\documentclass[12pt]{amsart}
\usepackage{amssymb, amsmath}

\newtheorem{theorem}{Theorem}[section]
\newtheorem{proposition}[theorem]{Proposition}
\newtheorem{corollary}[theorem]{Corollary}
\newtheorem{lemma}[theorem]{Lemma}
\newtheorem*{theoremnn}{Theorem}
\newtheorem*{corollarynn}{Corollary}

\theoremstyle{definition}
\newtheorem{definition}[theorem]{Definition}

\theoremstyle{remark}

\begin{document}

\newcommand{\fra}[1]{{\mathfrak{#1}}}
\newcommand{\Gr}{\text{Gr}}
\newcommand{\grh}{{\text{gr}}_{F}}
\newcommand{\ra}{\rightarrow}
\newcommand{\Ra}{\Rightarrow}
\newcommand{\surj}{\twoheadrightarrow}
\newcommand{\lra}{\longrightarrow}
\newcommand{\noi}{\noindent}
\newcommand{\PP}{\mathbf{P}}
\newcommand{\PPP}{\PP_{\text{sub}}}
\newcommand{\RR}{\mathbf{R}}
\newcommand{\NN}{\mathbf{N}}
\newcommand{\lef}{\mathbf{L}}
\newcommand{\hyp}{\mathbf{H}}
\newcommand{\ZZ}{\mathbf{Z}}
\newcommand{\CC}{\mathbf{C}}
\newcommand{\QQ}{\mathbf{Q}}
\newcommand{\bin}[2]{ {{#1} \choose {#2}} }
\newcommand{\OO}{\mathcal{O}}

\newcommand{\cH}{\mathcal{H}}
\newcommand{\MM}{\mathcal{M}}
\newcommand{\KK}{\mathcal{K}}
\newcommand{\II}{\mathcal{I}}
\newcommand{\LL}{\mathcal{L}}
\newcommand{\FF}{\mathcal{F}}
\newcommand{\GG}{\mathcal{G}}
\newcommand{\EE}{\mathcal{E}}
\newcommand{\cM}{\mathcal{M}}
\newcommand{\cN}{\mathcal{N}}
\newcommand{\cA}{\mathbf{A}}
\newcommand{\ti}[1]{\tilde{#1}}
\newcommand{\ef}{\rm{\ if\  }}
\newcommand{\eps}{\epsilon}
\newcommand{\sbl}{\vskip 3pt}
\newcommand{\lbl}{\vskip 6pt}
\newcommand{\rk} {\text{rank }}
\newcommand{\osect}{\mathbf{0}}
\newcommand{\HH}[3]{H^{{#1}} \big( {#2} , {#3} \big) }
\newcommand{\hh}[3]{h^{{#1}} \big( {#2} , {#3} \big) }
\newcommand{\OP}[1]{\OO_{\PP({#1})}}
\newcommand{\fall}{ \ \ \text{ for all } \ }
\newcommand{\rndup}[1]{ \ulcorner {#1} \urcorner }
\newcommand{\rndown}[1] {\llcorner {#1} \lrcorner}
\newcommand{\hgt}{\rm{height }}
\newcommand{\MI}[1]{\mathcal{J} ( {#1} ) }
\newcommand{\MIP}[1]{\mathcal{J}_+ ( {#1} ) }
\newcommand{\ZMI}[1]{\text{Zeroes}\big( \MI{ {#1} }
          \big)}
\newcommand{\BI}[1]{  \mathfrak{b} \big( {#1} \big) }
\newcommand{\ord}{\text{ord}}
\newcommand{\codim}{\text{codim}}
\newcommand{\mult}{\text{mult}}
\newcommand{\Supp}{\text{Supp}}
\newcommand{\defi}{{\text{def}}}
\newcommand{\pr}{\prime}
\newcommand{\QT}[2]{{#1}<{#2}>}
\newcommand{\Div}{\text{Div}}
\newcommand{\num}{ \equiv_{\text{num}} }
\newcommand{\lin}{\equiv}
\newcommand{\Qnum}{\equiv_{\text{num},\QQ}}
\newcommand{\Qlin}{\equiv_{\text{lin},\QQ}}
\newcommand{\dra}{\dashrightarrow}
\newcommand{\Bl}{\text{Bl}}
\newcommand{\lcm}{\text{lcm}}
\newcommand{\length}{\text{length}}
\newcommand{\lam}{\lambda}
\newcommand{\al}{\alpha}
\newcommand{\Spec}{\text{Spec}}
\newcommand{\HSp}{\text{HSp}}
\newcommand{\db}{\underline{\Omega}\,\dot{}}

\title{On Hodge spectrum and multiplier ideals}

\author{Nero Budur}
\address{Department of Mathematics \\University of
Illinois at Chicago \hfil\break\indent 851 South Morgan Street
(M/C 249)\\ Chicago, IL 60607-7045, USA} \email{nero@math.uic.edu}

\begin{abstract}
We describe a relation between two invariants that measure the
complexity of a hypersurface singularity. One is the Hodge
spectrum which is related to the monodromy and the Hodge
filtration on the cohomology of the Milnor fiber. The other is the
multiplier ideal, having to do with log resolutions.
\end{abstract}

\maketitle

\section{Introduction}
\label{intro}

Let $D$ be an effective divisor on a nonsingular complex variety
$X$. The  multiplier ideal $\MI{D}$ is a subsheaf of ideals of
$\OO_X$. It is related to the birational geometry of the pair
$(X,D)$ and measures in a subtle way the singularities of $D$. The
singularities of $D$ get "worse" if $\MI{D}$ is smaller (see
\cite{La}).

On the other hand, given a point $x\in X$ and a germ $f$ of a
hypersurface at $x$, singularity theory  brings into focus
standard topological invariants as the Milnor fiber and the
monodromy. An important analytic invariant is the canonical mixed
Hodge structure on the cohomology of the Milnor fiber (see
\cite{Ku}).

The natural question to ask is: How does the multiplier ideal
$\MI{D}$ reflect the standard constructions of singularity theory
of a germ of $D$ at a point $x$?

We will consider numerical invariants. The multiplier ideal
$\MI{D}$ comes with a finite increasing sequence of multiplier
ideals
$$\MI{D}=\MI{\al_0D}\subsetneq\MI{\al_1D}\subsetneq\ldots\subsetneq\MI{\al_nD}=\OO_X,$$
with $1=\al_0>\al_1>\ldots >\al_n >0$ and $\al_i\in\QQ$. The
$\al_i$ are called the jumping numbers of $(X,D)$ and provide some
measure of the complexity of the scheme defined by $\MI{D}$. The
jumping numbers seem to have appeared first in this form in
\cite{LV}. They are also related with the constants of
quasi-adjunction of \cite{Li}. The smallest jumping number is the
same as the log canonical threshold of $(X,D)$.  Let $x$ be a
point in $X$. Some jumping numbers $\al_i$ reflect closer than
others the difference at $x$ between $\MI{\al_iD}$ and
$\MI{\al_{i+1}D}$. We will call these  inner jumping numbers
(\ref{ij1}). To each $\al\in\QQ_{>0}$ we
 assign an inner jump multiplicity,
$$n_{\al,x}(D)\ge0,$$
with strict inequality iff $\al$ is an inner jumping number
(\ref{ij1}). If $D$ has only an isolated singularity at $x$, all
non-integral jumping numbers are inner jumping numbers.

The  Hodge spectrum of a local equation $f$ of $D$ at $x$  is a
fractional Laurent polynomial
$$\text{Sp}(f)=\sum_{\al\in\QQ}n_\al(f)t^\al,$$
with $n_\al(f)\in\ZZ$, and can be defined using the Hodge
filtration and the monodromy on the cohomology of the Milnor fiber
of $f$ (\ref{hs}). We call $n_{\al}(f)$ the multiplicity of $\al$
in the Hodge spectrum. There are various ways to encode the
multiplicities $n_\al(f)$ into a fractional Laurent polynomial.
Our definition of the spectrum coincides with that of M. Saito
\cite{Sa4}.

 The main result of this article is the
following.
\begin{theoremnn}
Let $X$ be a nonsingular quasi-projective variety of dimension
$m$. Let $D$ be an effective integral divisor on $X$ and $x\in D$
be a point. Let $f$ be any local equation of $D$ at $x$. Then for
any $\al\in (0,1]$, $n_{\al}(f)=n_{\al,x}(D).$
\end{theoremnn}

\begin{corollarynn} For all $\al\in (0,1]$,
\begin{enumerate}
\item $\al$ appears in the Hodge spectrum of $f$ if and only if $\al$ is an inner jumping
 number of $(X,D)$ at $x$;
\item (M. Saito \cite{Sa93}-0.7) the multiplicity $n_{\al}(f)$ is $\ge 0$;
\item (Varchenko \cite{Va}-\S 4; see also \cite{Ko}-\S 9) if $x$ is an isolated
singularity of $D$ and $\al\ne 1$, then, replacing $X$ by an open
neighborhood of $x$ if necessary, $\al$ appears in the Hodge
spectrum if and only if $\al$ is a jumping number.
\end{enumerate}
\end{corollarynn}

The method of proof of the main result is an application of the
theory of jet schemes and motivic integration. More specifically,
we make use of a deep theorem of J. Denef and F. Loeser
\cite{DL98} that shows how to compute the Hodge spectrum from a
log resolution of the pair $(X,D)$. In the first section, we
introduce the multiplier ideals and the inner jump multiplicities.
In the second section, we introduce the Hodge spectrum. In the
third section we review the result of \cite{DL98}. In the fourth
section, we prove the main theorem, up to a few lemmas about
quotient singularities which make the subject of the last section.

All varieties will be assumed to be defined over $\CC$.

{\it Acknowledgement.} I thank my advisor L. Ein for getting
me interested in this question and helping me with comments and
suggestions. I thank T. de Fernex, R. Lazarsfeld, M.
Musta\c{t}\v{a}, and K. Smith for several discussions and comments
on this work.

\section{The Multiplier Ideals}\label{multiplier ideals}
 We review the definition of multiplier ideals, jumping numbers and we introduce the inner jump multiplicities at a
 point.

 For a nonsingular variety $X$, let $\omega_X$
  denote the canonical sheaf. Let $K_X$ be the canonical divisor class,
that is
  $\omega_X=\OO_X(K_X)$. Let $\mu:Y\ra X$ be a proper birational morphism. The exceptional
  set of $\mu$ is denoted by $Ex(\mu)$. This is the set of points $\{y\in Y\}$ where $\mu$ is
  not biregular. We say that $\mu$ is a {
   resolution} if $Y$ is smooth. We say that a resolution $\mu$ is a {log resolution} of a
   finite collection of ideal sheaves $\mathfrak{a}_1,\ldots,\mathfrak{a}_s\subset\OO_X$ if we can write
    $\mu^{-1}\mathfrak{a}_i\cdot\OO_Y=\OO_Y(-H_i)$ where $H_i$ is an effective divisor on $Y$
    ($1\le i\le s$) such that $\cup_i H_i\cup Ex(\mu)$ is a divisor with
    simple normal crossings. Such a resolution always exists, by Hironaka. Let $K_{Y/X}=K_Y-\mu^*(K_X)$.
    If $D=\sum d_iD_i$ is an
    effective $\QQ$-divisor on $X$,
  where $D_i$ are the irreducible components of $D$ and $d_i\in\QQ$, the {round down}
  of $D$ is the integral divisor $\rndown{D}=\sum\rndown{d_i}D_i$. Here, $\rndown{d}$ is the
  biggest integer $\le d$.

  The following is well-known and follows from the Kawamata-Viehweg vanishing theorem
  (see \cite{La}, 9.2.18 and
  9.4.1):

\begin{theorem}\label{multiplier} Let $X$ be a nonsingular
quasi-projective variety, $\mathfrak{a}_1,\ldots,\mathfrak{a}_s$
$\subset\OO_X$ a finite collection of ideal sheaves, and
$d_1,\ldots,d_s\in\QQ_{>0}$ . Let $\mu:Y\ra X$ be a log resolution
of $(\mathfrak{a}_1,\ldots,\mathfrak{a}_s)$ and
$\mu^{-1}\mathfrak{a}_i\cdot\OO_Y=\OO_Y(-H_i)$. Then
 $\mu_*\OO_Y(K_{Y/X}-\rndown{d_1H_1+\ldots+d_sH_s})$ is independent of the choice of
 $\mu$ and
$R^i\mu_*\OO_Y(K_{Y/X}-\rndown{d_1H_1+\ldots+d_sH_s})=0$ for $i>0$.
\end{theorem}

\begin{definition} The {\it multiplier ideal} of  $\mathfrak{a}_1,\ldots,\mathfrak{a}_s\subset\OO_X$ and
$d_1,\ldots,d_s\in\QQ_{>0}$ is
$$\MI{\mathfrak{a}_1^{d_1}\cdot\ldots\cdot\mathfrak{a}_s^{d_s}}=\mu_*\OO_Y(K_{Y/X}-\rndown{d_1H_1+\ldots+d_sH_s})$$
for any log resolution $\mu$ of
$(\mathfrak{a}_1,\ldots,\mathfrak{a}_s)$.
\end{definition}

By Theorem \ref{multiplier}, the definition of
$\MI{\mathfrak{a}_1^{d_1}\cdot\ldots\cdot\mathfrak{a}_s^{d_s}}$ is
independent of the choice of the log resolution. The multiplier
ideals are indeed sheaves of ideals because
$\mu_*\omega_Y=\omega_X$. An important particular case is when
each $\mathfrak{a}_i=\OO_X(-D_i)$ with $D_i$ irreducible divisors
on $X$. The corresponding multiplier ideal will be denoted by
$\MI{D}$, where $D=\sum_{i=1}^s d_iD_i$. A log resolution of
$(\mathfrak{a}_1,\ldots,\mathfrak{a}_s)$ will be called a {\it log
resolution} of $(X,D)$. For $\al\in\QQ_{>0}$ define $\KK_\al
(D)=\MI{(1-\eps)\al\cdot D}/\MI{\al\cdot D},$ where $0<\eps<<1$.
\begin{definition} We say that $\al\in\QQ_{>0}$ is a {\it jumping number
of} $(X,D)$ if $\KK_\al(D)\ne 0$.
\end{definition}

Let $x\in D$ be a point. Next, we will see how some jumping
numbers $\al$ reflect closer the difference at $x$ between the
multiplier ideals $\MI{\al D}$ and $\MI{(1-\eps)\al D}$
($0<\eps\ll 1$). Denote by $\fra{a}$ the ideal sheaf $\OO_X(-D)$
and by $\fra{m}_x$ the ideal sheaf of the point $x$ in $X$. For
$\al\in\QQ_{>0}$, define
$\KK_{\al,x}(D)=\MI{\fra{a}^{(1-\eps)\al}}/\MI{\fra{a}^{(1-\eps)\al}\cdot\fra{m}_x^\delta},$
where $0<\eps\ll\delta\ll 1$. That this definition is independent
of choice of $\eps$ and $\delta$ will be clear from Proposition
\ref{ijn}. We will also see  that $\KK_{\al,x}(D)$ has support at
most $\{x\}$. Thus, we can make the following definition.

\begin{definition}\label{ij1} We say that $\al\in\QQ_{>0}$ is an {\it inner jumping number
of} $(X,D)$ at $x$ if $\KK_{\al,x}(D)\ne 0$. Let the {\it inner
 jump multiplicity} of $(X,D)$ at $x$ be $n_{\al,x}(D)= \dim_\CC\KK_{\al,x}(D)$.
\end{definition}

We fix the following notation for the rest of the section. Let
$\mu:Y\ra X$ be a log resolution of $(\fra{a},\mathfrak{m}_x)$.
Let $\mu^*(D)=\sum_{i}m_iE_i$, where $E_i$ are the irreducible
components. Let $J=\{i\ |\ m_i\ne 0\}$. For any subset $I\subset
J$, let $E^I=\bigcup_{i\in I}E_i$. Let $E^\emptyset=\emptyset$.
For a positive integer $d$, let $J_d=\{i\in J\ |\ d|m_i\}$. Let
${J_{d,x}}=\{i\in J_d\ |\ \mu(E_i)=\{x\}\ \}$. Fix from now
$\al\in\QQ_{>0}$. Write $\al=r/d$ with $r$ and $d$ nonnegative
integers such that $\gcd(r,d)=1$. Having fixed $\al$, define three
effective reduced divisors on $Y$,
\begin{align*}
E=E^{J_d}, F=E^{J_{d,x}}, G=E^{J_d\backslash {J_{d,x}}},
\end{align*}
such that $E=F+G$.

\begin{lemma}\label{doua_ideale} Let $H$ be an effective divisor
on $Y$ such that $\mu^{-1}\mathfrak{m}_x\cdot\OO_Y=\OO_Y(-H)$.
Then
\begin{align*}
\OO_Y(K_{Y/X}-\rndown{\al\mu^*D}+G)&=\OO_Y(K_{Y/X}-\rndown{(1-\eps)\al\mu^*D}-F)\\
&=\OO_Y(K_{Y/X}-\rndown{(1-\eps)\al\mu^*D+\delta H}),
\end{align*}
where $0<\eps\ll\delta\ll 1$.
\end{lemma}
\begin{proof}
The first equality is clear by the definition of $E=F+G$. Let
$k_i$ and $h_i$ be the coefficients of $E_i$ in $K_{Y/X}$ and,
respectively, $H$. Then the coefficient of $E_i$ in
$K_{Y/X}-\rndown{(1-\eps)\al\mu^*D+\delta H}$ is
$k_i-\rndown{(1-\eps)\al m_i+\delta h_i}$. The coefficient of
$E_i$ in $K_{Y/X}-\rndown{(1-\eps)\al\mu^*D}$ is
$k_i-\rndown{(1-\eps)\al m_i}$. If $E_i$ is not a component of
$E$, that is if $\al m_i\not\in\ZZ$, then for small enough $\eps$
and $\delta$ the two coefficients are the same. If $\mu(E_i)\ne x$
then $h_i=0$ and the two coefficients are again the same. Suppose
now that $\al m_i\in\ZZ$ and $\mu(E_i)=x$, that is $E_i$ is a
component of $F$. For small $\delta$ and $\eps\ll\delta$, for
example $\eps\le\delta h_i/(\al m_i)$, we have
$\rndown{(1-\eps)\al m_i+\delta h_i}=\al m_i=\rndown{(1-\eps)\al m_i}+1.$
Therefore
$k_i-\rndown{(1-\eps)\al m_i+\delta h_i}=k_i-\rndown{(1-\eps)\al
m_i}-1,$ so the two coefficients from above differ by 1, as
claimed. \qed\end{proof}
\begin{lemma}\label{dinamo}
\begin{enumerate}
\item $\mu_*\OO_Y(K_{Y/X}-\rndown{\al\mu^*D}+G)$ is independent of the log resolution and $R^i\mu_*\OO_Y(K_{Y/X}-\rndown{\al\mu^*D}+G)=0$ for $i>0$.
\item For $0<\eps\ll 1$, $\mu_*(\OO_F(K_{Y/X}-\rndown{(1-\eps)\al\mu^*D}))$ is independent of the log resolution and
$R^i\mu_*(\OO_F(K_{Y/X}-\rndown{(1-\eps)\al\mu^*D}))=0$ for $i>0$.
\item $\mu_*(\OO_G(K_{Y/X}-\rndown{\al\mu^*D}+G))$ is independent of the log resolution and $R^i\mu_*(\OO_G(K_{Y/X}-\rndown{\al\mu^*D}+G))=0$ for $i>0$.
\end{enumerate}
\end{lemma}
\begin{proof}
By Theorem \ref{multiplier} and Lemma \ref{doua_ideale},
$$\mu_*\OO_Y(K_{Y/X}-\rndown{\al\mu^*D}+G)=\MI{\mathfrak{a}^{(1-\eps)\al}\cdot\mathfrak{m}_x^\delta},$$
where $\mathfrak{a}=\OO_X(-D)$ and $0<\eps\ll\delta\ll 1$. This
shows (1). For all small $\eps>0$, consider the exact sequence
\begin{align*}
0\lra\OO_Y(A-F)\lra\OO_Y(A)\lra\OO_F(A)\lra 0,
\end{align*}
where $A=K_{Y/X}-\rndown{(1-\eps)\al\mu^*D}$. By (1) and Lemma
\ref{doua_ideale}, the first sheaf pushes forward to a multiplier
ideal on $X$ and has no higher direct images. The second sheaf
pushes forward to $\MI{(1-\eps)\al D}$ and has no higher direct
images by Theorem \ref{multiplier}. This shows (2). (3) is
similar. \qed\end{proof}

Remark that
$\KK_\al(D)=\mu_*(\OO_{E}(K_{Y/X}-\rndown{(1-\eps)\al\mu^*D})).$
Similarly, the proof of Lemma \ref{dinamo}-(2) shows that we
obtain the following interpretation of inner jumping numbers in
terms of log resolutions.

\begin{proposition}\label{ijn} With the above notation, we have:
\begin{enumerate}
\item $\KK_{\al,x}(D)=\mu_*(\OO_F(K_{Y/X}-
\rndown{(1-\eps)\al\mu^*D}))$, where $0<\eps\ll 1$.
\item
$n_{\al,x}(D)=\chi(Y,\OO_F(K_{Y/X}-\rndown{(1-\eps)\al\mu^*D})$,
where $\chi$ is the sheaf Euler characteristic.
\end{enumerate}
\end{proposition}

\begin{proposition}\label{inner} Let $\al\in\QQ_{>0}$. If $\al$ is an inner jumping number
of $(X,D)$ at $x$ then $\al$ is a jumping number of $(X,D)$.
\end{proposition}
\begin{proof} Consider the exact sequence
\begin{align*}
0\lra\OO_G(K_{Y/X}-\rndown{\al\mu^*D}+G)\lra\OO_E(K_{Y/X}-
\rndown{(1-\eps)\al\mu^*D})\\
\lra\OO_F(K_{Y/X}-\rndown{(1-\eps)\al\mu^*D})\lra 0.
\end{align*}
None of the three sheaves has higher direct images for $\mu$ and
the last two sheaves push-forward to $\KK_\al(D)$ and,
respectively, $\KK_{\al,x}(D)$. If $\KK_{\al,x}(D)\ne 0$ then
$\KK_\al(D)\ne 0$. \qed\end{proof}

\begin{proposition}\label{kalu} Let $X$ and $D$ be as above. If $x$ is an isolated singularity of $D$, and if $X$
is replaced by a small open neighborhood of $x$ if necessary, then
all non-integral jumping numbers are inner jumping numbers.
\end{proposition}
\begin{proof} Let $\al\in\QQ_{>0}$ but
$\not\in\ZZ$. We can assume that $F=E$. Then $\KK_{\al,x}(D)$ and
$\KK_\al(D)$ are the direct image under $\mu$ of the same sheaf. \qed\end{proof}

\section{The Hodge Spectrum}\label{hodge spectrum}
We recall the definition of the Hodge spectrum of a hypersurface
singularity. The Hodge spectrum was first defined in \cite{St77}
for an isolated singularity, and then in \cite{St89} for the
general case.

Let $\text{MHS}_\CC$ denote the abelian category of complex mixed
Hodge structures. For $n\in\ZZ$, let $\CC(n)$ denote the
1-dimensional complex mixed Hodge structure such that $\Gr^W_i=0$
for all $i\ne -2n$ and $\Gr_F^p=0$ for all $p\ne -n$. Here $W_.$
is the weight filtration and $F\,\dot{}$ is the Hodge filtration.

\begin{definition}
Let $H\in\text{MHS}_\CC$ and let $T$ be an automorphism of $H$ of
finite order. The {\it Hodge spectrum} of $(H,T)$ is the
fractional Laurent polynomial
$$\HSp (H,T)=\sum_{\al\in\QQ}n_\al\cdot t^\al,$$
with $n_\al=\dim\Gr ^{p}_FH_\lam$, where $p=\rndown{\al}$,
$H_\lam$ is the eigenspace of $T$ for the eigenvalue
$\lam=\exp(2\pi i\al)$, and $F\,\dot{}$ is the Hodge filtration.
\end{definition}
This definition extends to the Grothendieck group of the abelian
category of complex mixed Hodge structures with an automorphism
$T$ of finite order.

Let $f:(\CC^{m},0)\ra(\CC,0)$ be the germ of a non-zero
holomorphic function. Let  $M_t$ be the Milnor fiber of $f$
defined as
$$M_t=\{z\in\CC^{m}\ |\ |z|<\eps\ {\rm and\ }f(z)=t\}$$
for $0<|t|\ll\eps\ll 1$. The cohomology groups $H^*(M_t,\CC)$
carry a canonical mixed Hodge structure such that the semisimple
part $T_s$ of the monodromy acts as an automorphism of finite
order of these mixed Hodge structures (see \cite{St77} for $f$
with an isolated singularity only, \cite{Na} and \cite{Sa4} for
the general case). We will not need the actual construction of the
mixed Hodge structure on $H^*(M_t,\CC)$.

Let $$\HSp '(f)=\sum_{j\in\ZZ}(-1)^j\HSp(\ti{H}^{m-1+j}(M_t,\CC),
T_s),$$ where $\ti{H}$ stands for reduced cohomology.
\begin{definition}\label{hs} We call the {\it Hodge spectrum} of the germ $f$ the fractional
Laurent polynomial $\text{Sp}(f)$ obtained by replacing in $\HSp
'(f)$ the powers $t^\al$ with $t^{m-\al}$, for all $\al\in\QQ$. Denote by $n_{\al}(f)$ the coefficient of $t^\al$ in
$\text{Sp}(f)$.
\end{definition}

\section{The Hodge Spectrum after Denef and Loeser}\label{sec-2}
In \cite{DL98}, Denef and Loeser used the theory of motivic
integration to give an interpretation of the Hodge spectrum of a
germ of a hypersurface in terms of log resolutions.  The purpose
of this section is to review their result.

Let $\MM$ denote the pseudo-abelian category of Chow motives over
$\CC$ with complex coefficients (see \cite{Sc}). Let $K_0(\MM)$ be
the Grothendieck group of  $\MM$. On
$K_0(\MM)$ and $K_0(\text{MHS}_\CC)$, the tensor products induce
ring structures such that the Hodge realization $H:K_0(\MM)\ra
K_0(\text{MHS}_\CC)$ becomes a ring homomorphism. For $n$ an
integer, $H(\lef^n)=\CC(-n)$, where $\lef$ is the Lefschetz
motive. Let $G$ be a
finite abelian group and $\hat{G}$ its complex character group.
Define $\text{Sch}_G$ to be the category of varieties on which $G$
acts by automorphisms and such that the $G$-orbit of any closed
point is contained in an open affine subvariety.

\begin{theorem}\label{euler}(\cite{DL98}-1.3)  There exists a  map
$$\chi_c:\text{Sch}_G\times\hat{G}\lra K_0(\MM)$$
with the following properties:
\begin{enumerate}
\item $H\circ\chi_c(X,\al)=\sum_i(-1)^i[H^i_c(X)_\al]$, where
$H^*_c(X)_\al$ denotes subspace of the cohomology with compact
supports (with the Deligne-Hodge structure) on which $G$ acts by
multiplication with $\al$.
\item If $Y$ is a closed $G$-stable subvariety in $X\in\text{Sch}_G$, for all $\al\in\hat{G}$,
$\chi_c(X\backslash Y,\al)=\chi_c(X,\al)-\chi_c(Y,\al).$
\item If $X\in\text{Sch}_G$, and $U$ and $V$ are $G$-invariant open subvarieties of $X$, for all $\al\in\hat{G}$,
$\chi_c(U\cup V,\al)=\chi_c(U,\al)+\chi_c(V,\al)-\chi_c(U\cap V,\al).$
\item Suppose $X\in\text{Sch}_G$ and the $G$-action factors through a quotient $G\ra H$. If $\al$ is not in the image of $\hat{H}\ra\hat{G}$, then $\chi_c(X,\al)=0$.
\end{enumerate}
\end{theorem}

Let $X$ be a nonsingular and connected variety and let
$f:X\ra\mathbf{A}^1$ be a regular function on $X$. Let $D$ be the
effective integral divisor $\text{div}(f)$. Let $\mu:Y\ra X$ be a
log resolution of $(X,D)$. Let $\mu^*D=\sum_im_iE_i$, where $E_i$
are the irreducible components. Fix $\al\in [0,1)$ a rational
number, and write $\al=r/d$ where $r$ and $d$ are nonnegative
integers with $\gcd(r,d)=1$. If $\al=0$, let $d=1$. Let $J=\{i\ |\
m_i\ne 0\}$ and $J_d=\{i\in J\ |\ d|m_i\}$. For $I\subset J$, let
$E_I=\cap_{i\in I}E_i$ and $E^I=\cup_{i\in I} E_i$. Let
$E_\emptyset=Y$ and $E^\emptyset=\emptyset$. For $I\subset J$, let
$E_I^o=E_I-E^{J\backslash I}$.

Let $p:\ti{Y}\ra Y$ be the degree-$d$ cyclic cover obtained by
taking the $d$-th root of the local equations of the divisor
$R=\sum_{i\in J\backslash J_{d}}m_iE_i$ (see \cite{Ko}-\S 2 or
\cite{DL98}-3.2.3.). That is, locally $\ti{Y}$ is the
normalization of the subvariety of $Y\times\cA^1$ given by
$y^d=s$, where $s$ defines $R$. For any locally closed subset
$W\subset Y$ let $\ti{W}=p^{-1}(W)$.

The generator $1(\text{mod}\ d)$ of $G=\ZZ/d$ gives a canonical
isomorphism of the complex character group $\hat{G}$ with $\mu_d$,
the group of the $d$-th roots of unity. Denote the inverse of this
isomorphism by $\gamma$. For $Z\in\text{Sch}_G$, we deliberately
abuse the notation and write $\chi_c(Z,\al)$ in place of
$\chi_c(Z,\gamma(e^{2\pi i\al}))$.

\begin{definition}(\cite{DL98})\label{defs}
$S_{\al,x}=\sum_{I\subset J_d}\chi_c((E_I^o\cap \mu^{-1}(x))\ti{}\ ,\al)(1-\lef)^{|I|-1}.$
\end{definition}

Let $H:K_0(\MM)\ra K_0(\text{MHS}_\CC)$ be the Hodge realization.
\begin{theorem}\label{denef-loeser}(\cite{DL98}-4.3.1) Let $X$ be a nonsingular variety of dimension $m$. Let $D$ be an effective integral divisor on $X$ and $x\in D$ a point. Then for any local equation $f$ of $D$ at $x$,
$$(-1)^{m-1}+\HSp\;'(f)=(-1)^{m-1}\sum_{\al\in \QQ\cap [0,1)}\HSp(H(S_{\al,x}), {\rm Id})t^\al$$
\end{theorem}
 For $p\in\ZZ$, let $\grh^{p}:K_0(\text{MHS}_\CC)\ra\ZZ$
be the map such that to a complex mixed Hodge structure $H$
assigns $\grh ^p(H)=\dim\Gr^p_FH$, where $F\,\dot{}$ is the Hodge
filtration.
\begin{corollary}\label{gr(m-1)} For $\al\in \QQ\cap (0,1]$  the multiplicity of $\al$ in the Hodge spectrum
$\text{Sp}(f)$ of $f$ is
$n_{\al}(f)=(-1)^{m-1}\grh^{m-1}H(S_{1-\al,x}).$
\end{corollary}
\begin{proof} By definition, $n_{\al}(f)$ is the multiplicity in $\HSp '(f)$ of $m-\al$. From
Theorem \ref{denef-loeser} and the definition of $\HSp'(f)$, the
multiplicity of $m-\al$ in $\HSp'(f)$ is equal to
$(-1)^{m-1}\grh^{m-1}H(S_{1-\al,x})$, where $\grh^p$ is the Euler
characteristic described above. \qed\end{proof}

In working with $S_{\al,x}$ we will use the following result:
\begin{proposition}\label{cover} (\cite{Ko}-\S 2, \cite{DL98}-3.2.3.) Let $I\subset J_d$ and let $W=E_I$.
\begin{enumerate}
\item The $\ZZ/d$-action gives an eigensheaf decomposition
$${p}_*\OO_{\ti{W}}=\bigoplus_{0\le j<d}\OO_{W}\otimes\OO_Y(\rndown{\frac{j}{d}\mu^*D}),$$
and the $\ZZ/d$-action on each term is given by multiplication by
$e^{2\pi ij/d}$.
\item Above $W\cap E^{J\backslash J_d}$, the $\ZZ/d$-action factors through a $\ZZ/d'$-action for some $d'<d$.
\item $\ti{W}$ is locally for the Zariski topology quotient of a nonsingular
variety by a finite abelian group.
\end{enumerate}
\end{proposition}

\section{From Hodge Spectrum to Inner Jumping Numbers}\label{computation}
In this section we prove the Theorem from the introduction, up to
a few facts about quotient singularities which we leave for the
next section. Remark that the Corollary from the introduction
follows directly from the Theorem and Propositions \ref{inner} and
\ref{kalu}.

The following Mayer-Vietoris type of result can be proven by a
careful local computation.
\begin{lemma}\label{mv}
Let $X$ be a nonsingular variety of dimension $m$. Let
$E=\sum_{i\in I}E_i$ be a reduced divisor on $X$ with simple
normal crossings. For $L\subset I$, denote by $E_I$ the
intersection $\cap_{i\in L}E_i$. Then there exists an  exact
complex of coherent sheaves
\begin{equation}\tag{\ref{mv}.1}
0\lra\OO_E\lra\mathop{\bigoplus_{L\subset
I}}_{|L|=1}\OO_{E_I}\lra\ldots\lra\mathop{\bigoplus_{L\subset
I}}_{|L|=p}\OO_{E_I}\lra\ldots
\end{equation}
\end{lemma}

\begin{lemma}\label{omega}
With notation as above, assume  that $E_i$ have proper
support. Then, for a locally free sheaf $\FF$ on $X$,
$\chi(\omega_{E}\otimes\FF)=\sum_{\emptyset\ne L\subset I}\chi(\omega_{E_L}\otimes\FF),$
where  $\omega_.$ are dualizing sheaves and $\chi$ is the
sheaf Euler characteristic.
\end{lemma}
\begin{proof}
Tensor the exact complex (\ref{mv}.1) over $\OO_X$ with the dual
$\FF^{\vee}$ of $\FF$. Then, by additivity of the Euler
characteristic,
$$\chi(\OO_E\otimes\FF^{\vee})=\sum_{\emptyset\ne L\subset I}(-1)^{|L|-1}\chi(\OO_{E_L}\otimes\FF^{\vee} ).$$
The claim follows by Grothendieck-Serre duality. \qed\end{proof}

{\it Proof of Theorem.} Both $n_{\al}(f)$ and $n_{\al,x}(D)$ do
not change if $X$ is replaced by an open neighborhood $U$ of $x$.
Hence, we can assume that $f$ is a regular function on $X$ and
$D=\text{div}(f)$. Let $\mu:Y\ra X$ be a common resolution of the
point $x$ and of the divisor $D$. That is,
$\mu^{-1}\fra{m}_x\cdot\OO_Y=\OO_Y(-H)$ for some effective divisor
$H$ on $Y$, where $\fra{m}_x$ is ideal sheaf of $x$ in $X$, and
$Ex(\mu)+H+\text{div}(f\circ\mu)$ is a divisor with simple normal
crossings. Let $\mu^*D=\sum_i m_iE_i$. Let $J=\{\ i\ |\ m_i\ne
0\}$.

Let $\al$ be a rational number in $[0,1)$. Write $\al=r/d$, with
$r$ and $d$ nonnegative coprime integers. Let $J_d=\{\ i\in J\ |\
d|m_i\}$, and $J_{d,x}=\{\ i\in J_d\ |\ \mu(E_i)=x\}$. For
$I\subset J$, let $E_I=\cap_{i\in I}E_i$, and $E^I=\cup_{i\in
I}E_i$. Let $E_\emptyset=Y$ and $E^\emptyset=\emptyset$. Denote by
$E_I^o$ the complement $E_I-E^{J\backslash I}$. Recall that by the
previous section we have a degree-$d$ cyclic cover $p:\ti{Y}\ra
Y$, and for any locally closed subset $W\subset Y$ we defined
$\ti{W}$ as the inverse image of $W$ under $p$. Having fixed
$\al\in [0,1)$ a rational number, recall the definition of
$S_{\al,x}\in K_0(\MM)$ from \ref{defs}. Let $M=\{I\subset J_d\ |\ I\cap {J_{d,x}}\ne\emptyset\}$. Firstly, at motivic
level, we have
$$S_{\al,x}=\sum_{I\in M}\chi_c(\ti{E_I^o},\al)(1-\lef)^{|I|-1}=\star,$$
because $\mu^{-1}(x)$ is a divisor whose components have nonzero
multiplicity in $\mu^*D$. By
Theorem \ref{euler}-(2),
$$\star=\sum_{I\in M}[\chi_c(\ti{E_I},\al)-\chi_c((E_I\cap E^{J\backslash I})\ti{}\ ,\al)](1-\lef)^{|I|-1}.$$
By  Theorem \ref{euler}-(3), applied
to $\chi_c((E_I\cap E^{J\backslash I})\ti{}\ ,\al)$,
$$\star=\sum_{I\in M}\sum_{L\subset J\backslash I}(-1)^{|L|}\chi_c(\ti{E}_{I\cup L},\al)(1-\lef)^{|I|-1}.$$
If $L\not\subset J_d$ then $E_{I\cup L}$ is included in $E_I\cap
E^{J\backslash J_d}$, and by Proposition \ref{cover} the
$\ZZ/d$-action on $\ti{E}_{I\cup L}$ factors through a
$\ZZ/d'$-action with $d'<d$. Hence, by Theorem \ref{euler}-(4),
$\chi_c(\ti{E}_{I\cup L},\al)=0$. Therefore
$$\star=\sum_{I\in M}\sum_{L\subset J_d\backslash I}(-1)^{|L|}\chi_c(\ti{E}_{I\cup L},\al)(1-\lef)^{|I|-1}.$$
For an integer $p$, let $\grh^p:K_0(\text{MHS}_\CC)\ra\ZZ$ be the
map that sends a complex mixed Hodge structure $H$ to
$\dim_\CC\Gr_F^pH\in\ZZ$, where $F\,\dot{}$ is the Hodge
filtration. For $H\in K_0(\text{MHS}_\CC)$ and an integer $n$, one
has $\grh^p(H\cdot\CC(n))=\grh^{p+n}(H)$. Let
$\grh^{m-1}H(S_{\al,x})=\star\star$. Then
$$\star\star=\sum_{I\in M}\sum_{L
\subset J_d\backslash
I}\sum_{i=0}^{|I|-1}(-1)^{|L|+i}\binom{|I|-1}{i}\grh^{m-1-i}H(\chi_c(\ti{E}_{I\cup
L},\al).$$ By Proposition \ref{cover}, the varieties
$\ti{E}_{I\cup L}$ are locally for the Zariski topology quotients
of nonsingular varieties by finite abelian groups. Since $I\cap
{J_{d,x}}\ne \emptyset$, they are finite over closed subvarieties
of $\mu^{-1}(x)$ and therefore projective. Thus, by Proposition
\ref{qs}-(1), $\grh^{m-1-i}H(\chi_c(\ti{E}_{I\cup L},\al))=0$ if
$m-1-i>\dim \ti{E}_{I\cup L}$. Since the divisors $E_i$ are in
normal crossing, the dimension of $\ti{E}_{I\cup L}$ is
$m-|I|-|L|$. Hence in the summation above it suffices to take
$i=|I|-1$ and $L=\emptyset$, i.e. $\star\star=\sum_{I\in
M}(-1)^{|I|-1}\grh^{m-|I|}H(\chi_c(\ti{E}_{I},\al)).$ By
Proposition \ref{qs}-(2),
$$\star\star=\sum_{I\in M}(-1)^{|I|-1}\sum_i(-1)^i\dim
H^i(\ti{E}_I,\OO_{\ti{E}_I})_{1-\al},$$ where
$H^*(\ti{E}_I,\OO_{\ti{E}_I})_{1-\al}$ is the part of the
cohomology of $\OO_{\ti{E}_I}$ on which the canonical generator
$1\ (\text{mod}\ d)$ of $G=\ZZ/d$ acts by multiplication with
$e^{2\pi i(1-\al)}$. Since the covering $p:\ti{E}_I\ra E_I$ is a
finite morphism, we can replace $\OO_{\ti{E}_I}$ by
${p}_*\OO_{\ti{E}_I}$ in $\star\star$. By Proposition
\ref{cover}, the eigensheaf decomposition of ${p}_*\OO_{\ti{E}_I}$
gives
$$\star\star=\sum_{I\in M}(-1)^{|I|-1}\chi(\OO_{E_I}\otimes\OO_Y(\rndown{(1-\al)\mu^*D})),$$
where $\chi$ stands for the sheaf Euler characteristic. The
varieties $E_I$ are projective and  nonsingular by assumption. Let
$\omega_{E_I}$ be the canonical-dualizing sheaf. By duality,
\begin{align*}
\star\star&=(-1)^{m-1}\sum_{I\in M}\chi\left(\omega_{E_I}\otimes\OO_Y(-
\rndown{(1-\al)\mu^*D})\right)\\
&=(-1)^{m-1}\sum_{\emptyset\ne L\subset {J_{d,x}}}\left
[\chi\left(\omega_{E_L}\otimes\OO_Y(-
\rndown{(1-\al)\mu^*D})\right)+\right.\\
&\ \ \ \ \ \ \ \ +\sum_{\emptyset\ne K\subset J_d\backslash
{J_{d,x}}}\left.\chi (\omega_{E_{L\cup
K}}\otimes\OO_Y(-\rndown{(1-\al)\mu^*D})\right ].
\end{align*}
 Let $G=E^{J_d\backslash {J_{d,x}}}$. By Lemma \ref{omega} applied for each nonempty index set $L\subset {J_{d,x}}$ to the nonsingular projective varieties $E_L$ and the  reduced divisors $G\cap E_L$,
\begin{align*}
\star\star&=(-1)^{m-1}\sum_{\emptyset\ne L\subset {J_{d,x}}}\left [\chi\left(\omega_{E_L}\otimes\OO_Y(-\rndown{(1-\al)\mu^*D})\right)\right.+\\
&\left.\ \ \ \ \ \  +\chi\left(\omega_{G\cap
E_L}\otimes\OO_Y(-\rndown{(1-\al)\mu^*D})\right)\right],
\end{align*}
where $\omega_{G\cap E_L}$ is the dualizing sheaf of $G\cap E_L$.
By adjunction for $G\cap E_L$ in $E_L$,
\begin{align*}
\star\star&=(-1)^{m-1}\sum_{\emptyset\ne L\subset
{J_{d,x}}}\chi(\omega_{E_L}\otimes\OO_Y(-\rndown{(1-\al)\mu^*D}+G)).
\end{align*}
Let $F=E^{J_{d,x}}$. By Lemma \ref{omega},
$$\star\star=(-1)^{m-1}\chi(\omega_F\otimes\OO_Y(-\rndown{(1-\al)\mu^*D}+G)).$$
By adjunction for $F\subset Y$,
$\omega_F\simeq\OO_F\otimes\OO_Y(K_Y+F)$ and therefore
\begin{align*}
\star\star&=(-1)^{m-1}\chi(\OO_F\otimes\OO_Y(K_{Y/X}-\rndown{(1-\eps)(1-\al)\mu^*D})),
\end{align*}
for small $\eps>0$. We conclude by Proposition \ref{inner} that,
for $\al\in\QQ\cap [0,1)$,
$\grh^{m-1}H(S_{\al,x})=(-1)^{m-1}n_{1-\al,x}(D).$ In other words,
for $\al\in \QQ\cap (0,1]$,
$\grh^{m-1}H(S_{1-\al,x})=(-1)^{m-1}n_{\al,x}(D).$ By Corollary
\ref{gr(m-1)}, we obtained for $\al\in \QQ\cap [0,1)$ that
$n_{\al}(f)=n_{\al,x}(D).$

\section{Quotient Singularities}
In this section we prove the remaining facts needed in the proof
of the Theorem. These facts are conclusions we draw
 from the Hodge theory of varieties with at most quotient singularities.  We say that a variety $X$ has only
 {\it quotient singularities} if $X$ is locally for the Zariski topology quotient of a
nonsingular variety $Y$ by a finite group of automorphisms of $Y$.

\begin{theorem}\label{rational}(E. Viehweg) Let $X$ be a variety. Suppose $X$ has only quotient singularities. Then $X$ has rational singularities, that is, for every resolution $f:Y\ra X$, we have $f_*\OO_Y=\OO_X$ and $R^if_*\OO_Y=0$ for $i>0$.
\end{theorem}

Let $X$ be a variety and $G$ a finite group acting on $X$ by
automorphisms. Let $\FF$ be a sheaf of $\OO_X$-modules. We say
that $G$ acts on $\FF$, or that $\FF$ is a {\it $G$-sheaf}, if for
every $g\in G$ there exists a lifting $\lambda_g^\FF:\FF\ra
g^*\FF$ of $\OO_X\ra g^*\OO_X$, such that
$\lambda_1^\FF=\text{id}_\FF$ and
$\lambda_{hg}^\FF=g^*(\lambda_h^\FF)\circ\lambda_g^\FF$. $G$ acts
naturally on any locally free $\OO_X$-module. A morphism of
$G$-sheaves $\theta:\FF\ra\GG$ is said to be $G$-equivariant if it
commutes with the $G$-action, that is,
$g^*\theta\circ\lambda_g^\FF=\lambda_g^\GG\circ\theta$. The
$G$-sheaves with their $G$-equivariant morphisms form a category
which we denote by $\text{Mod}(X,G)$. If $f:Y\ra X$ is a
$G$-equivariant morphism, then the usual direct image functor
gives a functor $f_*:\text{Mod}(Y,G)\ra\text{Mod}(X,G)$.

Let $X$ be a variety with only quotient singularities and $G$ a
finite group acting on $X$ by automorphisms. Let $U$ be the open
subvariety consisting of the regular points of $X$, and denote by
$j:U\ra X$ the inclusion. Let $f:Y\ra X$ be a $G$-equivariant
resolution of $X$. Assume that $f$ gives an isomorphism of
$f^{-1}(U)$ with $U$. Denote by $\omega_Y$ and $\omega_U$ the
canonical invertible sheaves of $Y$ and, respectively, $U$.
\begin{theorem}\label{st}(\cite{St77}-1.11)  There is a $G$-equivariant isomorphism of
$G$-sheaves $ f_*\omega_Y\ra j_*\omega_U.$
\end{theorem}

For a variety $X$, let $D_{\text{filt}}(X)$ denote the derived
category of filtered complexes of $\OO_X$-modules with
differentials being differential operators of order $\le 1$
(\cite{DB81}-1.1). Let $D^b_{\text{filt},\text{coh}}(X)$ be the
full subcategory of $D_{\text{filt}}(X)$ of filtered complexes
$(K\,\dot{},F)$ such that, for all $i$, ${\Gr}_{F}^iK\,\dot{}$ is
an object of $D^b_{\text{coh}}(X)$, the derived category of
bounded complexes of $\OO_X$-modules with coherent cohomology
(\cite{DB81}-1.4). For an integer $p$ and a complex $K\,\dot{}$,
let $K[p]\,\dot{}$ denote the complex such that $K[p]^q=K^{p+q}$.
For our purpose, the Hodge filtration on the cohomology of the
variety $X$ is best understood via the Du Bois complex $\db_X$
constructed in \cite{DB81}. The following are the basic properties
of the Du Bois complex we need.

\begin{theorem}\label{dubois}(\cite{DB81}-3.21, 4.5, 4.12, 5.3, 5.9; \cite{GNPP}-V.3.5, V.4.1).
For every variety $X$ there exists a filtered complex
$(\underline{\Omega}\dot{}_X,F)$ in
$D^b_{\text{filt},\text{coh}}(X)$, canonical up to isomorphism,
with the following properties.
\begin{enumerate}
\item If $f:Y\ra X$ is a morphism of varieties, there exists in $D_{\text{filt}}(X)$ a natural morphism
$f^*:(\db_X,F)\ra Rf_*(\db_Y,F),$
where $Rf_*$ is the right derived functor of the filtered direct
image functor $f_*$.
\item The complex $\db_X$ is a resolution of the constant sheaf $\CC$ on $X$. If $X$ is
nonsingular then $\db_X$ is the usual De Rham complex of K\"ahler
differentials with the ``filtration b\^ete''.
\item If $X$ is proper,
$\hyp^{p+q}(X,\Gr_F^p\db_X)\simeq\Gr_F^pH^{p+q}(X,\CC).$
\item The complex $\Gr_F^p\db_X$ is zero for $p>\dim X$.
\item If $X$ has only quotient singularities, let $U$ be the open subvariety of regular
points of $X$ and $j:U\ra X$ the inclusion. Then  $j^*$ induces in
$D_{\text{filt}}(X)$ an isomorphism
$(\db_X,F)\ra j_*(\db_U,F).$
\item If a finite group $G$ acts on a variety $X$,
then $(\db_X,F)$ can be considered as an object in the
$G$-equivariant derived category $D_{\text{filt}}(X,G)$
(\cite{DB81}-5.4). In particular, $\db_X$ is a complex of
$G$-sheaves, the differential and the filtration commute with the
$G$-action. If $f:Y\ra X$ is a $G$-equivariant morphism, then
$f^*:(\db_X,F)\ra Rf_*(\db_Y,F)$ can be taken as a morphism in
$D_{\text{filt}}(X,G)$.
\end{enumerate}
\end{theorem}

 For a finite abelian group $G$, let $\hat{G}$ be the group of complex characters of $G$. If $H$ is
a $\CC$-vector space on which $G$ acts and $\al\in\hat{G}$ is a
complex character of $G$, let $H_\al=\{\ h\in H\ |\ gh=\al(g)h,\
\text{for all}\ g\in G\ \}.$ For $X\in\text{Sch}_G$ and $\al\in\hat{G}$, let
$\chi_c(X,\al)$ be the element in $K_0(\MM)$ from Theorem
\ref{euler}. Let
   $H:K_0(\MM)\ra K_0(\text{MHS}_\CC)$ be the Hodge realization.
   For an integer $p$, let $\grh^p:K_0(\text{MHS}_\CC)\ra\ZZ$ be
   the map that sends  $H\in \text{MHS}_\CC$ to $\dim_\CC\Gr^p_FH$, where $F\,\dot{}$
   is the Hodge filtration.

\begin{proposition}\label{qs}
Let $X\in\text{Sch}_G$ be a projective variety of dimension $n$.
For all $\al\in\hat{G}$, the following holds:
\begin{enumerate}
\item $\grh^pH(\chi_c(X,\al))=0$, if $p>n$.
\item If, in addition, $X$ has only quotient singularities,
$\grh^nH(\chi_c(X,\al))=\sum_i(-1)^i\dim_\CC H^i(X,\OO_X)_{\al^{-1}}.$
\end{enumerate}
\end{proposition}
\begin{proof} Let $H^*(X)=H^*(X,\CC)$.  The compatibility of
$\chi_c$ with the Hodge realization gives
$H(\chi_c(X,\al))=\sum_{i\in\ZZ}(-1)^i[H^i(X)_\al]$ in
$K_0(\text{MHS}_\CC).$ Hence,
$\grh^pH(\chi_c(X,\al))=\sum_{i\in\ZZ}(-1)^i\dim_\CC\Gr^p_FH^i(X)_\al,$
where $F$ is the Hodge filtration. By Theorem \ref{dubois}-(3) and
(6), we have a $G$-equivariant isomorphism $\Gr^p_FH^i(X)\simeq
\hyp^i(X,\Gr_F^p\db_X),$ where $\db_X$ is the filtered complex of
Du Bois. By Theorem \ref{dubois}-(4), if $p>n$ then
$\Gr_F^p\db_X=0$, hence the hypercohomology groups are also zero.
Thus $\Gr^p_FH^i(X)$ is zero, and in particular the subspace
$\Gr^p_FH^i(X)_\al$ is zero.
This proves (1).\\
\indent Suppose now that $X$ has only quotient singularities. Let
$U$ be the open subvariety of $X$ of regular points and $j:U\ra X$
the inclusion. By Theorem \ref{dubois}-(4),(5),(6), there is a
natural quasi-isomorphism, compatible with the $G$-action, of
$\Gr^n_F\db_X[n]$ and the complex with $j_*\omega_U$ concentrated
in degree zero, where $\omega_U$ is the canonical invertible sheaf
of $U$. Therefore, for any integer $i$, we have a natural
$G$-equivariant isomorphism $\hyp^i(X,\Gr^n_F\db_X[n])\simeq
H^i(X,j_*\omega_U).$ This gives a $G$-isomorphism
$\Gr^n_FH^i(X)\simeq H^{i-n}(X,j_*\omega_U).$ Let $f:Y\ra X$ be a
$G$-equivariant resolution of $X$. By Theorem \ref{st}, there is a
$G$-equivariant isomorphism $\Gr^n_FH^i(X)\simeq
H^{i-n}(X,f_*\omega_Y),$ where $\omega_Y$ is the canonical
invertible sheaf of $Y$. By the Grauert-Riemenschneider vanishing
theorem, $R^if_*\omega_Y=0$ for $i>0$ and so, by the Leray
spectral sequence, $H^{i-n}(X,f_*\omega_Y)=H^{i-n}(Y,\omega_Y)$.
Let $B$ denote the Serre duality pairing
$$H^{i-n}(Y,\omega_Y)\times H^{2n-i}(Y,\OO_Y)\lra H^n(Y,\omega_Y)=\CC.$$
Since $Y$ is nonsingular and projective, one knows that
$gB(.,.)=B(g.,g.)$, for all $g\in G$, and the action on
$H^n(Y,\omega_Y)$ is trivial. If $\al\in\hat{G}$ is a complex
character of $G$, then $\dim_\CC \Gr^n_FH^i(X)_\al=\dim_\CC
H^{2n-i}(Y,\OO_Y)_{\al^{-1}}.$  By Theorem \ref{rational},
$f_*\OO_Y=\OO_X$ and $R^if_*\OO_Y=0$, for $i>0$. Hence, there is a
$G$-equivariant isomorphism of $H^{2n-i}(Y,\OO_Y)$ and
$H^{2n-i}(X,\OO_X)$. Therefore, $\grh^nH(\chi_c(X,\al))=\sum_{i\in\ZZ}(-1)^i\dim_\CC
H^{2n-i}(X,\OO_X)_{\al^{-1}}$ which is equivalent with (2). \qed\end{proof}

\end{document}